\documentclass[11pt,reqno,oneside]{article}
\usepackage[latin1]{inputenc}
\usepackage{pstricks,pst-node,pst-tree}
\usepackage{amsmath}
\usepackage{amssymb}
\usepackage[lmargin=.55in, rmargin=.55in, tmargin=1in, bmargin=1in]{geometry}
\usepackage{xcolor}
\usepackage{graphics,graphicx}

\newtheorem{thm}{Theorem}
\newtheorem{cor}{Corollary}
\newtheorem{lmm}{Lemma}

\newtheorem{pro}{Proposition}

\definecolor{myblue}{rgb}{0,0.1,0.9}
\definecolor{myred}{rgb}{0.9,0.1,0}

\usepackage{ulem}
\def\pf{\noindent \emph{Proof.}\ }
\def\qed{{\quad\rule{1mm}{3mm}\,}}

\title{Distribution of external branch lengths in Yule trees}
\author{Filippo Disanto\thanks{Dipartimento di Matematica, Universit\`a di Pisa, Italy. Email: filippo.disanto@unipi.it} \quad and \quad Michael Fuchs\thanks{Department of Mathematical Sciences, National Chengchi University, Taipei 116, Taiwan. Corresponding author. Email: mfuchs@nctu.edu.tw.}}
\date{}

\begin{document}

\maketitle

\begin{abstract}
The Yule branching process is a classical model for the random generation of gene tree topologies in population genetics. It generates binary ranked trees---also called {\it histories}---with a finite number $n$ of leaves. 
We study the lengths $\ell_1 > \ell_2 > ... > \ell_k > ...$ of the external branches of a Yule generated random history of size $n$, where the length of an external branch is defined as the rank of its parent node.
When $n \rightarrow \infty$, we show that the random variable $\ell_k$, once rescaled as $\frac{n-\ell_k}{\sqrt{n/2}}$, follows a $\chi$-distribution with $2k$ degrees of freedom, with mean $\mathbb E(\ell_k) \sim n$ and variance $\mathbb V(\ell_k) \sim n \big(k-\frac{\pi k^2}{16^k} \binom{2k}{k}^2\big)$. Our results contribute to the study of the combinatorial features of Yule generated gene trees, in which external branches are associated with singleton mutations affecting individual gene copies. 
\end{abstract}

\vskip .2cm

%{\bf Mathematics Subject Classifications (2020)} ??? $\cdot$ ??? $\cdot$ ???

\section{Introduction}\label{sec:1}

The Yule distribution \cite{harding, yule} is a fundamental probability model of tree topologies, also called ``histories'', used in evolutionary analyses. Histories are full binary rooted trees, with a ranking of internal nodes that divides the tree in different layers (Fig. \ref{layers}A). The probabilistic features of Yule distributed histories have been subject of numerous investigations, with a particular interest on combinatorial properties that affect the frequency spectrum of mutations in population genetic tree models. A particular focus is on the length distribution of tree branches. Branch length can be seen as a discrete parameter---when only the number of tree layers spanned by a branch is considered---or as a time related quantity---when each tree layer is in turn considered with a length given by a continuous random variable. In the latter case, histories are called ``coalescent'' trees. While branch length of coalescent trees has been widely studied (see, e.g., \cite{blum,caliebe,dahmer,diehl,freund,fu,janson}), the discrete length of the edges of a random history has received less attention. 

In this paper, extending previous results \cite{DisantoAndWiehe}, we investigate the distribution of  the different lengths of the external branches---i.e., those branches ending with a leaf---of random histories of given size selected under the Yule model. External branch length is an important parameter to study as it relates to singleton mutations in the site frequency spectrum of population genetic trees. Denoting by $\ell_k$ the $k$th largest length of an external branch in a Yule distributed random history of $n$ leaves, our main finding is that, for every $k \geq 1$, the rescaled variable $\frac{n-\ell_k}{\sqrt{n/2}}$ follows asymptotically a $\chi$-distribution with $2k$ degrees of freedom, with convergence of all moments~(Theorem~\ref{teo}).

The paper is organized as follows. We introduce terminology and some useful properties of histories in Section~\ref{sec:2}, showing in particular that external branch lengths in random histories can also be analyzed in terms of peaks of random permutations. In Section \ref{long1}, we refine results of \cite{DisantoAndWiehe} finding a closed formula for the probability of the length, $\ell_1$, of the longest external branch in a random history of given size $n$ and a recurrence for calculating the probability of the $k$th largest length, $\ell_k$, of an external branch. For increasing $n$, the asymptotic distribution of the variables $\ell_1, \ell_2, ..., \ell_k, ...$ is finally examined in Section \ref{sec:4}. 

\section{Yule histories, external branches and non-peaks of permutations}\label{sec:2}

\begin{figure}[tpb]
\begin{center}
\begin{tabular}{c c c}
\includegraphics*[scale=0.68,trim=0 0 0 0]{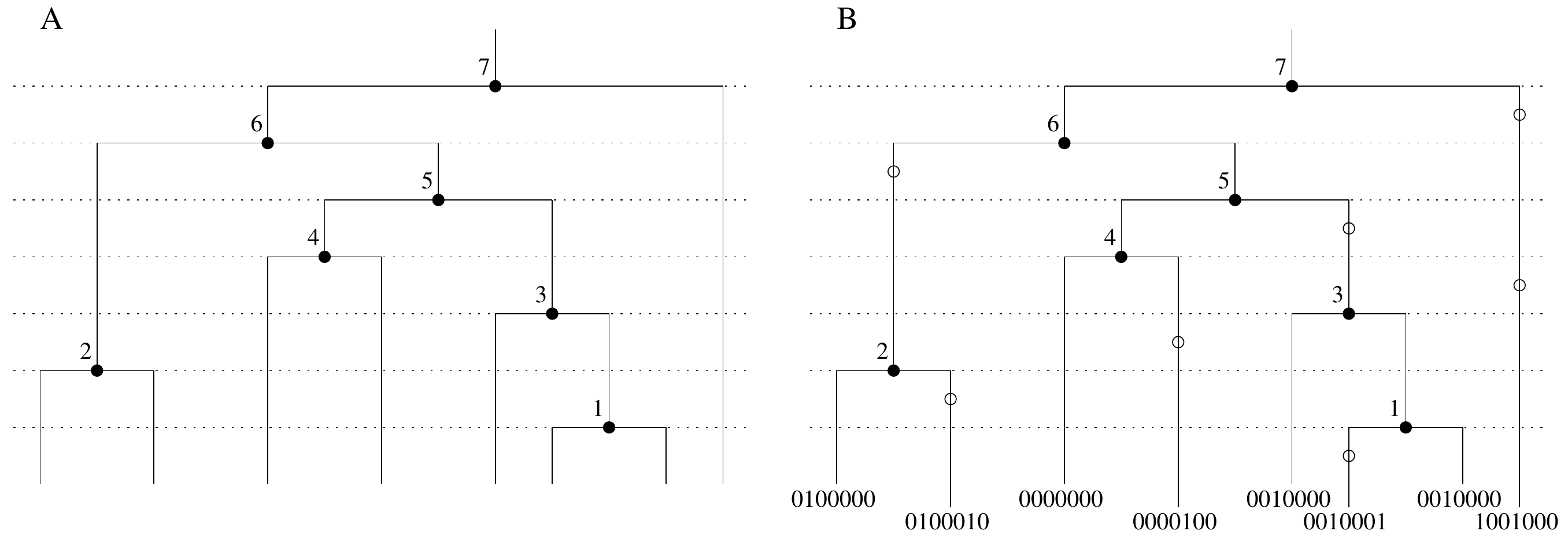}
\end{tabular}
\end{center}
\vspace{-.7cm}
\caption{{\small Histories and gene sequences. {\bf (A)} A history of size $n=8$.  The ranking of internal nodes decreases along any path going from the root to the leaves of the tree. The length of an external branch is the rank of its parent node. The different lengths of the external branches ordered by size are  $\ell_1 = 7 > \ell_2 = 4 > \ell_3 = 3 > \ell_4 = 2 > \ell_5 = 1$. 
{\bf (B)} The history depicted in A with leaves associated with genes represented as binary sequences with ancestral alleles of type $0$ and derived alleles of type $1$. A mutation (white circle) affects only the gene sequences associated with the leaves descending from the branch where it occurs. In this example, there is a mutation for each layer of the tree: the $i$th mutation (looking from top to bottom) changes the allele at the $i$th locus (position) of the gene.  
}} \label{layers}
\end{figure}

For a given positive integer $n$, a {\it history} \cite {rosenberg} of size $n$ is a full binary rooted tree with $n$ leaves and $n-1$ ranked internal nodes (Fig. \ref{layers}A). The rank of each internal node is defined by an integer label in $[1,n-1]$ bijectively associated with the node. The labeling decreases along any path from the root toward a leaf of the tree, determining a temporal ordering of the coalescent events---the merging of two edges---that characterize the branching structure of the tree. In a history of size $n$, there are $2n-1$ edges, or {\it branches}. A branch connecting an internal node and a leaf is said to be an {\it external} branch. The {\it length} of a branch is the difference between the rank of the nodes it connects. If the branch is external, then its length is simply the rank of its parent node.   

In Population Genetics, histories are tree structures that represent the evolution of individual genes from a common ancestor. Conditioning on a given history, an infinite sites model \cite{NielsenAndSlatkin} produces a set of mutations across the genes associated with the leaves of the tree. 
Roughly speaking, mutations occur as random events along the branches of the history (Fig. \ref{layers}B), with each branch containing a number of mutations that depends on its length, and with each mutation affecting only the set of gene copies descended from the branch it belongs to. 
In particular, a history with one or more ``long'' external branches will be associated with a biological scenario in which one or more gene copies will possess a ``large'' number of singleton mutations---i.e., mutations affecting only one individual. 
A random history of size $n$ selected under a proper null model distribution describes the evolutionary relationships of $n$ individual genes randomly sampled from a population under neutral evolution, and the length of the longest external branches in the random history relates to the largest number of singleton mutations that characterize single individuals in the sample.

In this paper, we focus on distributive properties of external branch length for random histories considered under a well known model of neutral evolution. More precisely, we will study external branch lengths ordered by size over random histories of size $n$ selected under the {\it Yule} probability model \cite{harding, yule}, or, equivalently, over random ordered histories of size $n$ selected uniformly at random. 
An {\it ordered} history of size $n$ is a plane embedding of a history of size $n$ in which subtrees carry a left-right orientation. The number of ordered histories of size $n$ is thus $(n-1)!$, and the Yule distribution over the set of histories of size $n$ is induced by the uniform distribution over the set of ordered histories of size $n$ by summing the probability $1/(n-1)!$ of each ordered history with the same underlying (un-ordered) history \cite{anconfig}. In particular, if $\text{c}(t)$ is the number of cherries (i.e., subtrees of size $2$) in a history $t$ of $n$ leaves, then $2^{n-1-\text{c}(t)}$ is the number of different plane embeddings of $t$, and therefore $2^{n-1-\text{c}(t)}/(n-1)!$ is the Yule probability of the history $t$ \cite{rosenberg}.

A series of combinatorial results on the lengths of external branches of uniformly distributed ordered histories (or Yule distributed histories) has been obtained in 
\cite{DisantoAndWiehe} in relationship with a study \cite{peaks} of the number of permutations of fixed size with a given set of peak entries, where the entry $\pi(i)$ is a {\it peak} in the permutation $( \pi(1), ..., \pi(i), ..., \pi(n) )$ when $i \neq 1, i \neq n$ and $\pi(i-1) < \pi(i) >\pi(i+1)$. Indeed, there exists a well known \cite{goulden} bijection that associates an ordered history $t$ of size $n$ with a permutation $\pi_t$ of the first $n-1$ positive integers. The mapping $t \rightarrow \pi_t$ can be described recursively by setting $\pi_t = (\pi_{t_L}, r(t), \pi_{t_R})$, where $r(t)$ is the (label of the) root of $t$, and $t_L, t_R$ are respectively the left and right root subtrees of $t$ (if any). In particular, ordering by size the different lengths $\ell_1 > \ell_2 > ... > \ell_k > ...$ of the external branches of $t$, the $k$th length, $\ell_k$, is easily seen to correspond to the $k$th largest non-peak entry in the permutation $\pi_t$.  
For example, if $t$ is the ordered history of size $n=8$ depicted in Fig. \ref{layers}, then $\pi_t = (2,6,4,5,3,1,7)$ has the following non-peak entries: $2, 4, 3, 1, 7$, which correspond to the different lengths $\ell_1 = 7 > \ell_2 = 4 > \ell_3 = 3 > \ell_4 = 2 > \ell_5 = 1$ of the external branches of $t$.
  By using the correspondence with non-peak entries of permutations, in the next section we calculate the probability of the varibale $\ell_k$ in an ordered history of size $n$ selected uniformly at random.    
  
\section{The probability of the $k$th external branch length}\label{long1}

Given an ordered history $t$ of size $n$, consider the different external branch lengths of $t$ ordered by size as $\ell_1 > \ell_2 > ... > \ell_k > ...$, where 
%$\lceil n/2 \rceil \leq k \leq n-1$ and 
$\ell_k \leq n-k$. 
%For a fixed $k$, if $t$ has sufficiently many leaves, then $\ell_k \geq 1+ \lceil (n-2k)/2  \rceil$.
%in order to obtain the minimum value of $\ell_k$, take for instance $t$ with one cherry whose parent node has rank $k$, with $k-1$ cherries whose parent nodes have rank larger than $k$, and with the remaining leaves arranged to form $\lfloor n/2 \rfloor - k$ cherries with parent nodes of rank smaller than $k$. 
As observed above, the value of $\ell_k$ corresponds to the $k$th largest non-peak entry in the associated permutation $\pi_t$. In this section, we study the number  $h_{n}(\ell_1=s_1, \ell_2=s_2, ..., \ell_k= s_k)$ of ordered histories of size $n$ in which $\ell_j = s_j$ for $j=1, ..., k$, which determines the probability $p_n(\ell_1=s_1, \ell_2=s_2, ..., \ell_k= s_k) = h_{n}(\ell_1=s_1, \ell_2=s_2, ..., \ell_k= s_k)/(n-1)!$. 

We start our calculations by using a result of \cite{peaks} for the number $\Pi_n(Q)$ of permutations of size $n \geq 3$ with peak entries matching the elements of a given set $Q \subseteq [3,n]$.
Fix $s_1, s_2, ..., s_{k-1}, s_k$ such that $n \geq s_1 > s_2 > ... > s_{k-1} > s_k$, and let $Z$ be a subset of the integers in the interval $[3,s_k-1]$. Then, by replacing $S=Z \cup [s_k+1,s_{k-1}-1] \cup [s_{k-1}+1,s_{k-2}-1] \cup ... \cup [s_2+1,s_1-1]$ and $K = n - s_1$ in Lemma 3.3 of \cite{peaks} (in which $k$ is the capital $K$ here), we find
\begin{eqnarray}\nonumber
&&\Pi_n(Z \cup [s_k+1,s_{k-1}-1] \cup [s_{k-1}+1,s_{k-2}-1] \cup ... \cup [s_2+1,s_1-1] \cup [s_1+1,n]) \\\nonumber 
&=& \Pi_n(S \cup [n-K+1,n]) 
= 2(K+1) \Pi_{n-1}(S \cup [n-K,n-1]) + K (K+1) \Pi_{n-2}(S \cup [n-K,n-2]) \\\nonumber
&=& 2(n-s_1+1) \Pi_{n-1}(Z \cup [s_k+1,s_{k-1}-1] \cup [s_{k-1}+1,s_{k-2}-1] \cup ... \cup [s_2+1,s_1-1] \cup [s_1,n-1]) \\\nonumber
&& + (n-s_1)(n-s_1+1) \Pi_{n-2}(Z \cup [s_k+1,s_{k-1}-1] \cup [s_{k-1}+1,s_{k-2}-1] \cup ... \cup [s_2+1,s_1-1] \cup [s_1,n-2]).
\end{eqnarray}
If we sum both sides of the latter equation over the possible subsets $Z$ of $[3,s_k-1]$, then we obtain
\begin{eqnarray}\label{sam}
&&\sum_Z \Pi_n(Z \cup [s_k+1,s_{k-1}-1] \cup [s_{k-1}+1,s_{k-2}-1] \cup ... \cup [s_2+1,s_1-1] \cup [s_1+1,n])  \\\nonumber
&=& 2(n-s_1+1) \, \sum_Z \Pi_{n-1}(Z \cup [s_k+1,s_{k-1}-1] \cup [s_{k-1}+1,s_{k-2}-1] \cup ... \cup [s_2+1,s_1-1] \cup [s_1,n-1]) \\\nonumber
&+& (n-s_1)(n-s_1+1) \, \sum_Z \Pi_{n-2}(Z \cup [s_k+1,s_{k-1}-1] \cup [s_{k-1}+1,s_{k-2}-1] \cup ... \cup [s_2+1,s_1-1] \cup [s_1,n-2]),
\end{eqnarray}
where the first sum 
%$\sum_Z \Pi_n(Z \cup [s_2+1,s_1-1] \cup [s_1+1,n])$ 
counts the permutations of size $n$ in which the first largest non-peak entry is $\ell_1=s_1$, the second largest non-peak entry is $\ell_2=s_2$, ..., and the $k$th largest non-peak entry is $\ell_k=s_k$. %non-peak entry satisfy $\ell_1=s_1, \ell_2=s_2$, ..., and $\ell_k=s_k$, while 
Similarly, the second and third sums 
%$\sum_Z\Pi_{n-1}(Z \cup [s_2+1,s_1-1] \cup [s_1,n-1])$ and $\sum_Z \Pi_{n-2}(Z \cup [s_2+1,s_1-1] \cup [s_1,n-2])$ 
count respectively the permutations of size $n-1$ and $n-2$ in which $\ell_1=s_2, \ell_2=s_3$, ..., and $\ell_{k-1}=s_{k}$. 
%the largest non-peak entry is $s_2$.
Note that when we set $k=1$ and $s_1=s$, we have $S=Z\subseteq [3,s-1]$ and the calculation above yields 
\begin{equation} \label{sami}
\sum_Z \Pi_n(Z \cup [s+1,n])  = 2(n-s+1) \, \sum_Z \Pi_{n-1}(Z \cup [s,n-1]) + (n-s)(n-s+1) \, \sum_Z \Pi_{n-2}(Z \cup [s,n-2]),
\end{equation}
where the first sum 
%$\sum_S\Pi_n(S \cup [s+1,n])$ 
counts the permutations of size $n$ in which the largest non-peak entry is $\ell_1=s$, while the second and third sums 
%$\sum_S\Pi_{n-1}(S \cup [s,n-1])$ and $\sum_S \Pi_{n-2}(S \cup [s,n-2])$ 
count respectively the permutations of size $n-1$ and $n-2$ in which the largest non-peak entry is strictly smaller than $s$, that is,  $\ell_1 < s$. By rewriting (\ref{sam}) and (\ref{sami}) in terms of ordered histories, we find 
\begin{eqnarray}\nonumber
h_{n+1}(\ell_1=s_1, \ell_2=s_2, ..., \ell_k= s_k) &=& 2(n-s_1+1) \, h_{n}(\ell_1=s_2, \ell_2=s_3, ..., \ell_{k-1}=s_{k}) \\\label{joint1}
&& + (n-s_1)(n-s_1+1) \, h_{n-1}(\ell_1=s_2, \ell_2=s_3, ..., \ell_{k-1}=s_{k}) 
\end{eqnarray}
and 
\begin{equation}\label{sd}
h_{n+1}(\ell_1=s) = 2(n-s+1) \, h_n(\ell_1 < s) + (n-s)(n-s+1) \, h_{n-1}(\ell_1<s),
\end{equation}
where $h_i(\ell_1 < s) \equiv \sum_{j < s} h_i(\ell_1 = j)$.

Because $h_{n+1}(\ell_1=s) = h_{n+1}(\ell_1 < s+1) - h_{n+1}(\ell_1 < s)$, Eq. (\ref{sd}) yields the recurrence
$h_{n+1}(\ell_1 < s+1) = h_{n+1}(\ell_1<s) +  2(n-s+1) \, h_n(\ell_1 < s) + (n-s)(n-s+1) \, h_{n-1}(\ell_1<s),$
which, by replacing $n+1$ by $n$ and $s+1$ by $s$, reads as 
\begin{equation}\label{ress}
h_{n}(\ell_1 < s) = h_{n}(\ell_1<s-1) +  2(n-s+1) \, h_{n-1}(\ell_1 < s-1) + (n-s)(n-s+1) \, h_{n-2}(\ell_1<s-1),
\end{equation} 
where $h_n(\ell_1<s) = 0$ if $s = \lceil n/2 \rceil$ ($\ell_1$ is at least $\lceil n/2 \rceil$), and $h_n(\ell_1<s) = (n-1)!$ if $s = n$ ($\ell_1$ is at most $n-1$).
In particular, when $\lceil n/2 \rceil \leq s \leq n \geq 3$, we have
\begin{equation}\label{p4}
h_n(\ell_1 < s) 
= \frac{(s-1)! \, (s-2)! \, (2s-n) \, (2s-n-1)}{(2s-n)!}
\end{equation} 
as the right-hand side---say $r(n,s)$---of the latter equation satisfies the same recurrence  (\ref{ress}) given for $h_{n}(\ell_1 < s)$. Indeed, $r(n,\lceil n/2 \rceil) = 0$ and $r(n,n) = (n-1)!$. Furthermore, assuming $\lceil n/2 \rceil < s < n$, a simple calculation shows that $r(n,s) = r(n,s-1) + 2(n-s+1) \, r(n-1,s-1) + (n-s)(n-s+1) \, r(n-2,s-1)$, where we note that all the factorials in $r(n,s-1), r(n-1,s-1)$, and $r(n-2,s-1)$ are well defined being of the form $m!$ with $m \geq 0$.

The next proposition summarizes our enumerative results from a probability point of view.
\begin{pro}\label{firstprop}
Let $n\geq 3$. If $p_n(\ell_1=s)$ denotes the probability of $\ell_1=s$ in an ordered history of size $n$ selected uniformly at random, then
\begin{equation}\label{prop1a}
p_n(\ell_1=s) %= \frac{h_n(\ell_1<s+1)-h_n(\ell_1<s)}{(n-1)!}
= \frac{(s-1)! (s-2)! (4 n s + s - n^2 - n -3 s^2)}{(2s-n)! \, (n-1)!},
\end{equation}
where  $\lceil n/2 \rceil \leq s \leq n-1$. Furthermore, the joint probability $p_n(\ell_1=s_1, \ell_2=s_2,..., \ell_k=s_k)$ of $\ell_1=s_1, \ell_2=s_2$, ..., and $\ell_k = s_k$ in an ordered history of size $n$ selected uniformly at random satisfies the recurrence
\begin{eqnarray} \label{joint2}
p_n(\ell_1=s_1, \ell_2=s_2,..., \ell_k=s_k) &=& \frac{2(n - s_1)}{n-1} p_{n-1}(\ell_1=s_2, \ell_2=s_3,..., \ell_{k-1}=s_k) \\\nonumber 
&& + \frac{(n - s_1)(n - s_1 - 1)}{(n-1)(n-2)} p_{n-2}(\ell_1=s_2, \ell_2=s_3,..., \ell_{k-1}=s_k),
\end{eqnarray}
with initial condition given by (\ref{prop1a}).
\end{pro}
\pf  Equation (\ref{prop1a}) follows from (\ref{p4}) as  $p_n(\ell_1=s) = [h_n(\ell_1<s+1)-h_n(\ell_1<s)]/(n-1)!$. The recurrence in (\ref{joint2}) is obtained  by replacing $n+1$ by $n$ in (\ref{joint1}) and dividing both sides of the resulting equation by $(n-1)!$. \qed

By summing over the possible values of $\ell_1,...,\ell_{k-1}$ the joint probability $p_n(\ell_1=s_1, \ell_2=s_2,..., \ell_k=s_k)$ yields for $k \geq 2$ the probability of $\ell_k = s_k$ in random ordered history of $n$ leaves:
\begin{equation}\label{summa}
p_n(\ell_k=s_k) = \sum_{s_1=s_k+k-1}^{n-1} \sum_{s_2=s_k+k-2}^{s_1-1} ... \sum_{s_i=s_k+k-i}^{s_{i-1}-1} ... \sum_{s_{k-1}=s_k+1}^{s_{k-2}-1} p_n(\ell_1=s_1, \ell_2=s_2,..., \ell_k=s_k).
\end{equation}
For instance, if $k=2$, then we obtain
\begin{eqnarray}\label{prop1b}
p_n(\ell_2=s_2) &=& \sum_{s_1=s_2+1}^{n-1} p_n(\ell_1=s_1, \ell_2=s_2) \\\nonumber
&=&\sum_{s_1=s_2+1}^{n-1}  \frac{2(n - s_1)}{n-1} p_{n-1}(\ell_1=s_2) + \frac{(n - s_1)(n - s_1 - 1)}{(n-1)(n-2)} p_{n-2}(\ell_1=s_2) \\\nonumber
&=& \frac{2 p_{n-1}(\ell_1=s_2)}{n-1} \sum_{s_1=s_2+1}^{n-1}  (n - s_1) + \frac{p_{n-2}(\ell_1=s_2)}{(n-1)(n-2)} \sum_{s_1=s_2+1}^{n-1} (n - s_1)(n - s_1 - 1), 
%\\\nonumber
%&=& \frac{p_{n-1}(\ell_1=s_2)}{n-1} (s_2 + s_2^2 - n - 2 s_2 n + n^2) + \frac{p_{n-2}(\ell_1=s_2)}{3(n-1)(n-2)} (-s_2^3 + 3 s_2^2 (-1 + n) + s_2 (-2 + 6 n - 3 n^2) + n (2 - 3 n + n^2)) 
\end{eqnarray}
which can be used together with (\ref{prop1a}), when $n\geq5$ and $s_2$ is in the range $\lceil n/2 \rceil-1\leq s_2 \leq n-2$. Similarly, if $k=3$, then we have 
\begin{eqnarray}\label{prop1c}
p_n(\ell_3=s_3) &=& \sum_{s_1=s_3+2}^{n-1} \sum_{s_2=s_3+1}^{s_1-1} p_n(\ell_1=s_1, \ell_2=s_2, \ell_3=s_3) \\\nonumber
&=& \sum_{s_1=s_3+2}^{n-1} \sum_{s_2=s_3+1}^{s_1-1} \frac{2(n - s_1)}{n-1} p_{n-1}(\ell_1=s_2, \ell_2=s_3) + \frac{(n - s_1)(n - s_1 - 1)}{(n-1)(n-2)} p_{n-2}(\ell_1=s_2, \ell_2=s_3) \\\nonumber
%= \sum_{s_1=s_3+2}^{n-1} \sum_{s_2=s_3+1}^{s_1-1} \frac{2(n - s_1)}{n-1} \bigg[\frac{2(n - 1 - s_2)}{n-2} p_{n-2}(\ell_1=s_3) + \frac{(n - 1 - s_2)(n - s_2 - 2)}{(n-2)(n-3)} p_{n-3}(\ell_1=s_3)\bigg]  \\\nonumber
% + \sum_{s_1=s_3+2}^{n-1} \sum_{s_2=s_3+1}^{s_1-1} 
%\frac{(n - s_1)(n - s_1 - 1)}{(n-1)(n-2)} \bigg[\frac{2(n -2 - s_2)}{n-3} p_{n-3}(\ell_1=s_3) + \frac{(n -2 - s_2)(n - s_2 - 3)}{(n-3)(n-4)} p_{n-4}(\ell_1=s_3) \bigg] \\\nonumber
&=&
\frac{4p_{n-2}(\ell_1=s_3)}{(n-1)(n-2)} \sum_{s_1=s_3+2}^{n-1} \sum_{s_2=s_3+1}^{s_1-1} (n-s_1)(n-1-s_2) 
\\\nonumber
&&+ \frac{2 p_{n-3}(\ell_1=s_3)}{(n-1)(n-2)(n-3)} \sum_{s_1=s_3+2}^{n-1} \sum_{s_2=s_3+1}^{s_1-1} 
(n-s_1)(n-s_2-2)(2n-2-s_2-s_1) \\\nonumber
&&+ \frac{p_{n-4}(\ell_1=s_3)}{(n-1)(n-2)(n-3)(n-4)}  \sum_{s_1=s_3+2}^{n-1} \sum_{s_2=s_3+1}^{s_1-1} (n-s_1)(n-s_1-1)(n-2-s_2)(n-s_2-3),
\end{eqnarray}
which can be coupled with (\ref{prop1a}), when $n\geq7$ and $\lceil n/2 \rceil-2\leq s_3 \leq n-3$.

\section{Asymptotic distribution of the $k$th external branch length} \label{sec:4}
In this section, we derive distributive properties of the random variable $\ell_k$---the $k$th largest external branch length---considered over ordered histories of size $n$ selected under the uniform distribution. We start by considering the case $k=1$, and then generalize to arbitrary values of $k$.

By dividing Eq. (\ref{p4}) by the number $(n-1)!$ of ordered histories of size $n$, we obtain the probability
\[
p_n(\ell_1<s)=\frac{(s-1)!(s-2)!}{(2s-n-2)!(n-1)!},\, \lceil n/2\rceil<s\leq n,
\]
or alternatively, with $u=s-1$,
\begin{equation}\label{dis-func}
p_n(\ell_1\leq u)=\frac{u!(u-1)!}{(2u-n)!(n-1)!},\, \lceil n/2\rceil\leq u\leq n-1.
\end{equation}
Our first result is the following local limit theorem.

\begin{lmm} \label{xzc}
When $n\rightarrow \infty$,
\begin{itemize}
\item[(a)] the probability $p_n(\ell_1=\lfloor n-x\sqrt{n/2}\rfloor)$ admits an asymptotic expansion of the form 
\[
p_n(\ell_1=\lfloor n-x\sqrt{n/2}\rfloor)=\frac{x}{\sqrt{n/2}}e^{-x^2/2}(1+o(1))+{\mathcal O}\left(\frac{e^{-x^2/2}}{n}\right)
\]
uniformly for $0\leq x\leq x^{*}\equiv n^{1/7}$.
\item[(b)] Furthermore,
\[
p_n(\ell_1\leq n-x^{*}\sqrt{n/2})={\mathcal O}\left(e^{-n^{2/7}/2}\right),
\]
with $x^{*}$ as defined in part (a).
\end{itemize}
\end{lmm}
\pf For part (a), first assume that $x\leq x^{*}$ is such that $u\equiv n-x\sqrt{n/2}$ is a non-negative integer smaller than $n$. Then,
%, for $n$ sufficiently large, $u-1 \geq \lceil n/2\rceil$ and 
Eq. (\ref{dis-func}) yields
\[
p_n(\ell_1=u)=p_n(\ell_1 \leq u)-p_n(\ell_1 \leq u-1)=\frac{u!(u-1)!}{(2u-n)!(n-1)!}-\frac{(u-1)!(u-2)!}{(2u-2-n)!(n-1)!}.
\]
Plugging in Stirling's formula $z! \sim z^z e^{-z} \sqrt{2 \pi z} \left(1 + \frac{1}{12 z} + \frac{1}{288 z^2} - \frac{139}{51840 z^3} - ... \right)$ gives the (complete) asymptotic expansion
\[
p_n(\ell_1=u)\sim \frac{x}{\sqrt{n/2}}e^{-x^2/2}\left(1+\sum_{d=1}^{\infty}\frac{q_d(x)}{n^{d/2}}\right),
\]
where $q_{d}(x)$ is a polynomial of degree $3d$. Thus, for the given range of $x$, $q_{d}(x)={\mathcal O}(n^{3d/7})$ and consequently
\[
\frac{q_{d}(x)}{n^{d/2}}={\mathcal O}(n^{3d/7-d/2})=o(1).
\]
This shows that $\sum_{d=1}^{k}\frac{q_d(x)}{n^{d/2}}=o(1)$ for every choice of $k$ and the claimed expansion (without the last term) holds for this case. Note that the case $u=n$, i.e., $x=0$, is trivially covered as $p_n(\ell_1=n)=0$.

Next, if $u$ is not an integer, then $\lfloor u\rfloor=u+{\mathcal O}(1)=n-x\sqrt{n/2}+{\mathcal O}(1)=n-(x+{\mathcal O}(1/\sqrt{n}))\sqrt{n/2}$, and thus we are in the first case with $x$ replaced by $\tilde{x}=x+{\mathcal O}(1/\sqrt{n})$. Hence,
\begin{eqnarray}\nonumber
p_n(\ell_1=\lfloor u\rfloor)
&=& \frac{\tilde{x}}{\sqrt{n/2}}e^{-\tilde{x}^2/2}(1+o(1))
= \frac{x+{\mathcal O}(1/\sqrt{n})}{\sqrt{n/2}}e^{-x^2/2+o(1)}(1+o(1)) \\\nonumber
&=& \frac{x}{\sqrt{n/2}}e^{-x^2/2}(1+o(1))+{\mathcal O}\left(\frac{e^{-x^2/2}}{n}\right),
\end{eqnarray}
which establishes the claim also in this case.

For part (b), we are interested in $p_n(\ell_1\leq \lfloor n-x^{*}\sqrt{n/2} \rfloor)$. Starting from (\ref{dis-func}), we use  Stirling's approximation $\log(z) = z \log(z)-z+(1/2) \log(2 \pi z) + o(1)$ to expand $\log(p_n(\ell_1 \leq u)) = \log(u!)+\log((u-1)!)-\log((2u-n)!)-\log((n-1)!)$ as 
\begin{small}
$$\frac{1}{2} (2 (n-2 u) \log (2 u-n)-\log (2  u -  n)-2 n \log (n-1)+\log (n-1)+(2 u-1) \log (u-1)+2 u \log (u)+\log (u)) + o(1).$$
\end{small}
Then, we plug in $u = \lfloor n-x^{*}\sqrt{n/2} \rfloor = n-n^{1/7}\sqrt{n/2}-c_n$, where $c_n$ is the fractional part of $ n-n^{1/7}\sqrt{n/2}$, and replace the resulting terms of the form $\log(n+f(n))$ by $\log(n)+f(n)/n-f(n)^2/n^2$ (where $f(n)/n \rightarrow 0$). Simple algebraic manipulations finally give
$$
\log(p_n(\ell_1 \leq \lfloor n-x^{*}\sqrt{n/2} \rfloor)) = -\frac{n^{2/7}}{2} + o(1),$$
which shows the claim. \qed

From the previous lemma, we obtain the following proposition that describes the asymptotic distribution of the random variable $\ell_1$ considered over ordered histories of size $n$ selected uniformly at random.
\begin{pro}\label{limit-law}
As $n\rightarrow\infty$,
\[
\frac{n-\ell_1}{\sqrt{n/2}}\stackrel{d}{\longrightarrow}{\rm Rayleigh}(1)
\]
with convergence of all moments. In particular, the mean and the variance of $\ell_1$ satisfy respectively
\begin{equation}\label{mean-var}
{\mathbb E}(\ell_1)\sim n \qquad \text{and} \qquad{\mathbb V}(\ell_1) \sim \left(1-\frac{\pi}{4} \right) n.
\end{equation}
\end{pro}

\pf Fix an $x\geq 0$. In order to prove the limit law, we have to show that, when $n\rightarrow \infty$, the probability of $(n-\ell_1)/\sqrt{n/2}\leq x$ converges to $1-e^{-x^2/2}$, which is the cumulative function of the Rayleigh distribution with parameter $1$. We first write  
\begin{align}
p_n\left(\frac{n-\ell_1}{\sqrt{n/2}}\leq x\right)&=p_n(n-x\sqrt{n/2}\leq \ell_1)=p_n(\lceil n-x\sqrt{n/2}\rceil\leq \ell_1)=
\sum_{s=\lceil n-x\sqrt{n/2}\rceil}^{n}p_n(\ell_1=s) \\\label{calc2}
&=\sum_{t=0}^{\tilde{x}}p_n(\ell_1= n-t\sqrt{n/2}),
\end{align}
where the latter sum is in steps of size $\sqrt{2/n}$ and $\tilde{x}=x+{\mathcal O}(1/\sqrt{n})$ is such that $n-\tilde{x}\sqrt{n/2}=\lceil n-x\sqrt{n/2}\rceil$.
For $n$ sufficiently large, we can assume $\tilde{x} \leq x \leq n^{1/7}$ and thus use part (a) of the lemma writing (\ref{calc2}) as
\begin{equation}\label{calc}
\sum_{t=0}^{\tilde{x}} \frac{t}{\sqrt{n/2}}e^{-t^2/2}(1+o(1))+{\mathcal O}\left(\frac{e^{-t^2/2}}{n}\right) = \sum_{t=0}^{\tilde{x}} \frac{t}{\sqrt{n/2}}e^{-t^2/2}(1+o(1)) + \sum_{t=0}^{\tilde{x}} {\mathcal O}\left(\frac{e^{-t^2/2}}{n}\right).
\end{equation}
Because the $1+o(1)$ factor in the second sum of (\ref{calc}) holds uniformly, it can be put in front of the sum obtaining
$$\sum_{t=0}^{\tilde{x}} \frac{t}{\sqrt{n/2}}e^{-t^2/2}(1+o(1)) = (1+o(1)) \sum_{t=0}^{\tilde{x}} \frac{t}{\sqrt{n/2}}e^{-t^2/2} = (1+o(1)) \sum_{t=0}^{x} \frac{t}{\sqrt{n/2}}e^{-t^2/2} + o(1),$$
where the upper limit in the last sum is now $x$. Moreover, the third sum in (\ref{calc}) can be bounded as
$$\sum_{t=0}^{\tilde{x}} {\mathcal O}\left(\frac{e^{-t^2/2}}{n}\right) = {\mathcal O}\left(\sum_{t=0}^{\infty} \frac{e^{-t^2/2}}{n}\right) = o(1).$$
%{\mathcal O}\left(\frac{1}{\sqrt{n}}\sum_{t=0}^{\infty} \frac{e^{-t^2/2}}{\sqrt{n/2}}\right)$$
Hence, for $n\rightarrow \infty$, the probability $p_n\left(\frac{n-\ell_1}{\sqrt{n/2}}\leq x\right)$ converges to the Riemann sum  $\sum_{t=0}^{x} \frac{t}{\sqrt{n/2}}e^{-t^2/2}$ with step size $dt=\sqrt{2/n}$, which can be approximated by the integral  $\int_{0}^{x}t e^{-t^2/2}{\rm d}t=1-e^{-x^2/2},$
as claimed. 

By a similar approach, one can also show that all moments converge. Starting from
\[
{\mathbb E}\left(\frac{n-\ell_1}{\sqrt{n/2}}\right)^m=\sum_{s=0}^{n}\left(\frac{n-s}{\sqrt{n/2}}\right)^m p_n(\ell_1=s),
\]
%\[
%{\mathbb E}\left(\frac{n-\ell_1}{\sqrt{n/2}}\right)^m=\sum_{s=\lceil n/2\rceil}^{n-1}\left(\frac{n-s}{\sqrt{n/2}}\right)^m p_n(\ell_1=s),
%\]
we replace $s$ by $s=n-x\sqrt{n/2}$ and break the sum into two parts obtaining
$$\sum_{x=0}^{\sqrt{2n}} x^m p_n(\ell_1=n-x\sqrt{n/2} )=\sum_{0 \leq x< n^{1/7}} x^m p_n(\ell_1=n-x\sqrt{n/2} ) +\sum_{n^{1/7} \leq x \leq \sqrt{2n}} x^m p_n(\ell_1=n-x\sqrt{n/2} ) \equiv \Sigma_1+\Sigma_2,$$
where all the sums proceed in steps of size $\sqrt{2/n}$. 
%and the upper limit $u$ is such that $n-u\sqrt{n/2} = \lceil n/2 \rceil$. 
For $\Sigma_2$, by part (b) of the lemma, we have
\[
\Sigma_2={\mathcal O}\left(n^{m/2}e^{-n^{2/7}/2}\right)=o(1).
\]
For $\Sigma_1$, by part (a) of the lemma, we have
\[
\Sigma_1=(1+o(1))\sum_{0\leq x<n^{1/7}}\frac{x^{m+1}}{\sqrt{n/2}}e^{-x^2/2}+{\mathcal O}\left(n^{-1}\sum_{0\leq x<n^{1/7}}e^{-x^2/2}\right).
\]
Here, the Riemann sum in $\Sigma_1$ can be approximated by the integral $\int_{0}^{n^{1/7}}x^{m+1}e^{-x^2/2}{\rm d}x$, which converges to $\int_{0}^{\infty}x^{m+1}e^{-x^2/2}{\rm d}x$. Overall,
\[
{\mathbb E}\left(\frac{n-\ell_1}{\sqrt{n/2}}\right)^m\stackrel{n \to \infty}{\longrightarrow} \int_{0}^{\infty}x^{m+1}e^{-x^2/2}{\rm d}x
\]
which proves the claimed convergence of moments. Finally, (\ref{mean-var}) follows from this convergence by straightforward computation.\qed

Note that when the limit distribution is uniquely determined by its moment sequence (which is the case for the Rayleigh distribution), convergence of all moments implies weak convergence. Although the second part of the proof of the latter proposition suffices to show that also the first claim holds true, we decided to provide the calculations for the convergence in distribution with the aim of improving the readability of the remaining part of the proof.

In the following, our goal is to show that, for an arbitrary fixed value of $k \geq 1$, the random variable $\ell_k$ follows asymptotically a $\chi$ distribution with $2k$ degrees of freedom. Indeed, note that the Rayleigh distribution found for the case $k=1$ is a $\chi$ distribution with $2$ parameters.

The next lemma describes the solution to the recurrence (\ref{joint2}) for the joint probability $p_n(\ell_1=s_1, \ell_2=s_2,..., \ell_k=s_k)$ and a formula for the probability $p_n(\ell_k=s_k)$ given in (\ref{summa}) in terms of the probability of $\ell_1 = s_k$ in trees of size smaller than or equal to $n$.
\begin{lmm}\label{muevu}
By setting 
$\mu_n(x)\equiv\frac{2x}{n-1}$ and $\nu_n(x)\equiv\frac{x(x-1)}{(n-1)(n-2)}$, we have 
\begin{equation}\label{piox}
p_n(\ell_1=s_1, \ell_2=s_2,..., \ell_k=s_k) = 
\sum_{\omega}\left(\prod_{\ell=0}^{k-2}\omega^{[\ell]}_{n-n_{\omega,\ell}-\ell}\left(n-n_{\omega,\ell}-\ell-s_{\ell+1}\right)\right)p_{n-n_{\omega,k-1}-k+1}(\ell_1 = s_k),
\end{equation}
where the sum runs over all words $\omega = \omega^{[0]}\cdots\omega^{[k-2]}$ of length $k-1$ with letters from the alphabet $\{\mu,\nu\}$, and $n_{\omega,\ell}$ is the number of $\nu$ in the first $\ell$ letters of $\omega$  (with $n_{\omega,0}=0$). With the same notation, we also have 
\begin{equation} \label{multi-sum}
p_n(\ell_k=s_k) = \sum_{s_1=1}^{s_k^{*}}\sum_{s_2=s_1}^{s_k^*}\cdots\sum_{s_{k-1}=s_{k-2}}^{s_k^{*}}\sum_{\omega}\left(\prod_{\ell=0}^{k-2}\omega^{[\ell]}_{n-n_{\omega,\ell}-\ell}(s_{\ell+1}-n_{\omega,\ell})\right)
p_{n-n_{\omega,k-1}-k+1}(\ell_1=s_k),
\end{equation}
where $s_k^{*}\equiv n-k+1-s_k$.
\end{lmm} 
\begin{figure}[tpb]
\begin{center}
\includegraphics*[scale=0.63,trim=0 0 0 0]{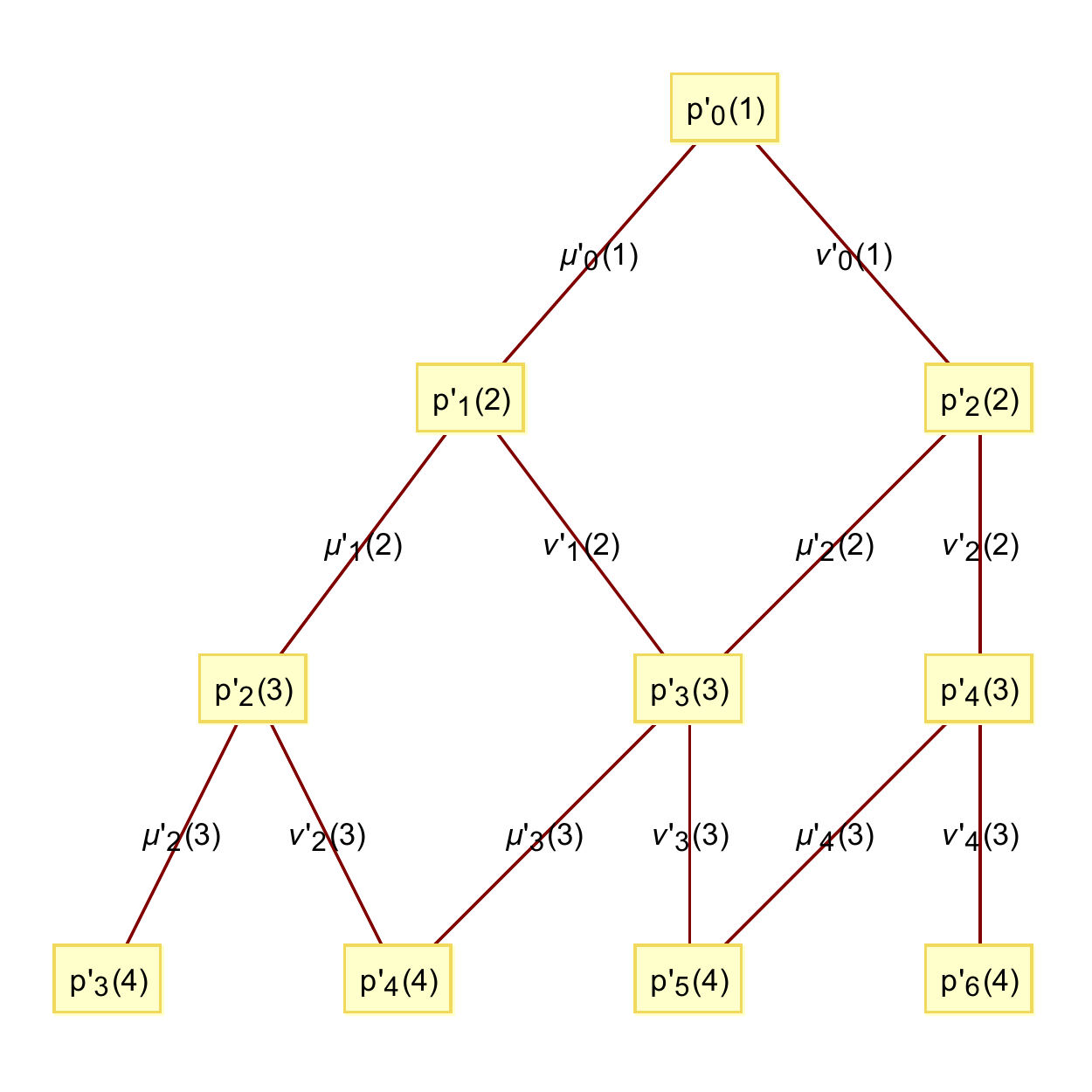}
\end{center}
\vspace{-.7cm}
\caption{{\small Schematic diagram of the first three iterative steps of the procedure (\ref{proce}) for calculating $p'_0(1)=p_n(\ell_1=s_1, \ell_2=s_2,..., \ell_k=s_k)$.
}} \label{alb}
\end{figure} 
\pf For a fixed $n$ and $k$, set $p'_i(j) \equiv p_{n-i}(\ell_1=s_j,...,\ell_{k-j+1}=s_k)$, $\mu'_i(j) \equiv \frac{2(n-i-s_{j})}{n-i-1}$, and $\nu'_i(j) \equiv \frac{(n-i-s_{j})(n-i-1-s_j)}{(n-i-1)(n-i-2)}$. The recurrence (\ref{joint2}) finds $p'_0(1)=p_n(\ell_1=s_1, \ell_2=s_2,..., \ell_k=s_k)$ by iteratively computing  
\begin{equation}\label{proce}
p'_i(j)=\mu'_i(j) \, p'_{i+1}(j+1) + \nu'_i(j) \, p'_{i+2}(j+1).
\end{equation}
The procedure ends after $k-1$ steps, that is, when we obtain terms of the form $p_{n-x}(\ell_1=s_k)=p'_x(k)$, for a certain value of $x$. For $k=4$, the diagram in Fig. \ref{alb} depicts the three iterations needed for evaluating $p'_0(1)$. The latter quantity is calculated as the sum of the probabilities at the bottom of the diagram, each multiplied by the sum of the words of length $k-1$ over the alphabet  $\{\mu',\nu'\}$ that encode the different paths connecting the corresponding leaf node to the root of the diagram.
More precisely, for arbitrary values of $n$ and $k$, we have
$$p'_0(1) = \sum_{\omega}\left(\prod_{\ell=0}^{k-2}\omega^{[\ell]}_{n_{\omega,\ell}+\ell}\left(\ell+1\right)\right)p'_{n_{\omega,k-1}+k-1}(k),$$
where the sum runs over all words $\omega = \omega^{[0]}\cdots\omega^{[k-2]}$ of length $k-1$ with letters from the alphabet $\{\mu',\nu'\}$, and $n_{\omega,\ell}$ is the number of $\nu'$ in the first $\ell$ letters of $\omega$  (with $n_{\omega,0}=0$). By replacing indices, the latter formula is equivalent to that claimed in (\ref{piox}).

Finally, plugging (\ref{piox}) into (\ref{summa}) yields 
\begin{align*}
p_n(\ell_k=s_k)=&\sum_{s_1=s_k+k-1}^{n-1}\sum_{s_2=s_{k}+k-2}^{s_1-1}\cdots\sum_{s_{k-1}=s_k+1}^{s_{k-2}-1}\sum_{\omega}\left(\prod_{\ell=0}^{k-2}\omega^{[\ell]}_{n-n_{\omega,\ell}-\ell}
\left(n-n_{\omega,\ell}-\ell-s_{\ell+1}\right)\right)\nonumber\\
&\hspace*{7cm}\times p_{n-n_{\omega,k-1}-k+1}(\ell_1=s_k).
\end{align*}
By setting $s_{\ell}^{*}=n-\ell+1-s_{\ell}$ for $\ell=1,...,k$, 
the right-hand side can be written as 
$$\sum_{s_1^{*}=1}^{s_k^{*}}\sum_{s_2^{*}=s_1^{*}}^{s_k^*}\cdots\sum_{s_{k-1}^{*}=s_{k-2}^{*}}^{s_k^{*}}\sum_{\omega}\left(\prod_{\ell=0}^{k-2}\omega^{[\ell]}_{n-n_{\omega,\ell}-\ell}(s_{\ell+1}^{*}-n_{\omega,\ell})\right)
p_{n-n_{\omega,k-1}-k+1}(\ell_1=s_k),
$$
which gives (\ref{multi-sum}).
\qed

With the same notation used above, we now provide two more useful lemmas.

\begin{lmm}\label{ll-longest}
For $s_k=\lfloor n-x\sqrt{n/2}\rfloor$, we have
\[
\sum_{s_1=1}^{s_k^{*}}\sum_{s_2=s_1}^{s_k^*}\cdots\sum_{s_{k-1}=s_{k-2}}^{s_k^{*}}\prod_{\ell=0}^{k-2}\mu_{n-\ell}(s_{\ell+1})=\frac{x^{2k-2}}{2^{k-1}(k-1)!}+{\mathcal O}\left(\frac{1+ x^{2k-3}}{\sqrt{n}}\right)
\]
uniformly for $0 \leq x \leq \sqrt{2n}$.
\end{lmm}
\pf Note that  
$$\sum_{s_1=1}^{s_k^{*}}\sum_{s_2=s_1}^{s_k^*}\cdots\sum_{s_{k-1}=s_{k-2}}^{s_k^{*}}\prod_{\ell=0}^{k-2}\mu_{n-\ell}(s_{\ell+1})=\frac{2^{k-1}   \sum_{s_1=1}^{s_k^{*}} s_1 \sum_{s_2=s_1}^{s_k^*} s_2 \cdots\sum_{s_{k-1}=s_{k-2}}^{s_k^{*}} s_{k-1}}{(n-1) \cdots (n-k+1)}= \frac{2^{k-1}r(s_k^{*})}{n^{k-1}}+{\mathcal O}\left(\frac{r(s_k^{*})}{n^{k}}\right),$$
where $r(z)$ is the polynomial $r(z)\equiv \sum_{s_1=1}^{z}s_1 \sum_{s_2=s_1}^{z} s_2 \cdots\sum_{s_{k-1}=s_{k-2}}^{z} s_{k-1}$.
In order to determine the asymptotic behavior of $r(z)$, we rely on Faulhaber's formula:
\begin{equation}\label{Faul}
\sum_{m=1}^{N}m^t=\frac{1}{t+1}\sum_{k=0}^{t}\binom{t+1}{k}B_k(N+1)^{t+1-k}\stackrel{N \rightarrow \infty}{\sim} \frac{N^{t+1}}{t+1}\stackrel{N \rightarrow \infty}{\sim}\int_{1}^{N}x^t{\rm d}x,
\end{equation}
where $B_k$ denotes the $k$-th Bernoulli number. In particular, we use the fact that, for a given polynomial $p(u)=\alpha_k u^k + ... + \alpha_1 u + \alpha_0$, the polynomial $\sum_{u=a}^b p(u) = \sum_{u=1}^b p(u) - \sum_{u=1}^{a-1} p(u)$ has its term $\frac{\alpha_k b^{k+1}}{k+1}$ with the highest power in $b$ and its term $-\frac{\alpha_k a^{k+1}}{k+1}$ with the highest power in $a$ matching those that appear in the integral $\int_{a}^b p(z) {\rm d}z$. As a consequence, if we substitute each sum in $r(z)$ by an integral sign, we then find a polynomial $\int_{1}^{z} s_1 {\rm d}s_{1} \int_{s_1}^{z} s_2 {\rm d}s_{2} \cdots\int_{s_{k-2}}^{z} s_{k-1} \, {\rm d}s_{k-1}$ with the same leading term of $r(z)$. 
Furthermore, by a simple induction on $k$ one can show that  
$\int_{z_{k+1}}^z z_{k} {\rm d} z_{k} \cdots \int_{z_3}^z z_2 {\rm d} z_2 \int_{z_2}^z z_1 {\rm d} z_1 = \frac{1}{2^{k}} \sum_{i=0}^{k} \frac{(-1)^{i} z^{2 k - 2 i} z_{k+1}^{2 i}}{i! (k-i)!},$
and therefore the leading term of $r(z)$ is that of $\frac{1}{2^{k-1}}\sum_{i=0}^{k-1} \frac{(-1)^{i} z^{2 k - 2 - 2 i}}{i! (k-1-i)!},$ that is, $\frac{x^{2k-2}}{2^{k-1} (k-1)!}$. Hence,
$$\sum_{s_1=1}^{s_k^{*}}\sum_{s_2=s_1}^{s_k^*}\cdots\sum_{s_{k-1}=s_{k-2}}^{s_k^{*}}\prod_{\ell=0}^{k-2}\mu_{n-\ell}(s_{\ell+1})=\frac{(s_k^{*})^{2k-2}}{n^{k-1}(k-1)!}+{\mathcal O}\left(\frac{(s_k^{*})^{2k-3}}{n^{k-1}}+\frac{r(s_k^{*})}{n^{k}}\right)$$
By plugging $s_k^{*}=x\sqrt{n/2}+{\mathcal O}(1)$ into the latter asymptotic formula and performing a straightforward expansion, we obtain the claimed result. \qed

The next result shows that Lemma \ref{ll-longest} gives the main term of the multiple sum in (\ref{multi-sum}).
\begin{lmm}\label{le4}
For $s_k=\lfloor n-x\sqrt{n/2}\rfloor$, we have
\[
\sum_{s_1=1}^{s_k^{*}}\sum_{s_2=s_1}^{s_k^*}\cdots\sum_{s_{k-1}=s_{k-2}}^{s_k^{*}}\prod_{\ell=0}^{k-2}\omega^{[\ell]}_{n-n_{\omega,\ell}-\ell}(s_{\ell+1}-n_{\omega,\ell})={\mathcal O}\left(\frac{1+ x^{2k-1}}{\sqrt{n}}\right)
\]
uniformly for $0 \leq x \leq \sqrt{2n}$ and for all words $\omega = \omega^{[0]}\cdots\omega^{[k-2]}$ of length $k-1$ with letters from the alphabet $\{\mu,\nu\}$ different from the word whose letters are all equal to $\mu$.
\end{lmm}
\pf Assume that $\omega$ has $m\geq 1$ letters equal to $\nu$. Then, since $\nu_n(x)$ is a quadratic polynomial, by again using Faulhaber's formula (\ref{Faul}), we obtain that
\[
\sum_{s_1=1}^{s_k^{*}}\sum_{s_2=s_1}^{s_k^*}\cdots\sum_{s_{k-1}=s_{k-2}}^{s_k^{*}}\prod_{\ell=0}^{k-2}\omega^{[\ell]}_{n-n_{\omega,\ell}-\ell}(s_{\ell+1}-n_{\omega,\ell})
=\frac{r(s_{k}^{*})}{q(n)},
\]
where $r(z)$ is a polynomial of degree $m+2k-2$ and $q(z)$ is a polynomial of degree $m+k-1$. Thus, by setting $s_k^{*}=x\sqrt{n/2}+{\mathcal O}(1)$, we obtain that
\[
\sum_{s_1=1}^{s_k^{*}}\sum_{s_2=s_1}^{s_k^*}\cdots\sum_{s_{k-1}=s_{k-2}}^{s_k^{*}}\prod_{\ell=0}^{k-2}\omega^{[\ell]}_{n-n_{\omega,\ell}-\ell}(s_{\ell+1}-n_{\omega,\ell})
=\frac{r(s_{k}^{*})}{q(n)}={\mathcal O}\left(\frac{1+ x^{m+2k-2}}{n^{m/2}}\right).
\]
From this the result follows by observing that $x \leq \sqrt{2n}$.\qed

\bigskip

From the last three lemmas, we can now deduce the following generalization of Lemma \ref{xzc}.
\begin{cor} \label{coro1}
When $n \rightarrow \infty$,
\begin{itemize}
\item[(a)] the probability $p_n(\ell_k = \lfloor n-x\sqrt{n/2}\rfloor)$ admits an asymptotic expansion of the form 
\[
p_n(\ell_{k}=\lfloor n-x\sqrt{n/2}\rfloor)=\frac{x^{2k-1}}{2^{k-1}(k-1)!\sqrt{n/2}}e^{-x^2/2}(1+o(1))+{\mathcal O}\left(\frac{e^{-x^2/2}}{n}\right)
\]
uniformly for $0\leq x\leq x^{*}\equiv n^{1/7}$.
\item[(b)] Furthermore,
\[
p_n(\ell_{k}\leq n-x^{*}\sqrt{n/2})={\mathcal O}\left(n^{k-1}e^{-n^{2/7}/2}\right)
\]
with $x^{*}$ as defined in part (a).
\end{itemize}
\end{cor}

\pf First, note that for any given word $\omega$ of length $k-1$ over the alphabet $\{\mu,\nu\}$ (in the sense of Lemma \ref{muevu}), we have 
\[
p_{n-n_{\omega}-k+1}(\ell_1=\lfloor n-x\sqrt{n/2}\rfloor)=p_{n-n_{\omega}-k+1}(\ell_1=\lfloor n-n_{\omega}-k+1-\tilde{x}\sqrt{(n-n_{\omega}-k+1)/2}\rfloor),
\] 
where $\tilde{x}=x+{\mathcal O}(1/\sqrt{n})$. As a consequence, by applying part (a) of Lemma 1 with $x$ replaced by $\tilde{x}$ and $n$ replaced by $n-n_{\omega}-k+1$, it follows that part (a) of Lemma~\ref{xzc} also holds when $p_n$ is replaced by $p_{n-n_{\omega}-k+1}$.
Moreover, also part (b) of Lemma~\ref{xzc} holds true when $p_n$ is replaced by $p_{n-n_{\omega}-k+1}$. Indeed, from (\ref{dis-func}), we find
\begin{align*}
p_{n-n_{\omega}-k+1}(\ell_1\leq n^{*})&= \frac{n^{*}! (n^{*}-1)!}{(2 n^{*}- n+n_{\omega}+k-1  )! (n-n_{\omega}-k)!}\\ &=\frac{(n-1)\cdots(n-n_{\omega}-k+1)}{(2n^{*}-n+n_{\omega}+k-1)\cdots(2n^{*}-n+1)} \cdot \frac{n^{*}!(n^{*}-1)!}{(2n^{*}-n)!(n-1)!} = {\mathcal O}(p_n(\ell_1 \leq n^{*})),
\end{align*}
where $n_{\omega}\equiv n_{\omega,k-1}$ and $n^{*} \equiv \lfloor n-x^{*}\sqrt{n/2}\rfloor$. 
%{\bf(To see that part (a) of Lemma~\ref{xzc} holds when $p_n$ is replaced by $p_{n-n_{\omega}-k+1}$, are we taking the first display of this proof and apply part (a) of Lemma 1 with $x$ replaced by $\tilde{x}$ and $n$ replaced by $n-n_{\omega}-k+1$? Or it is better to say that, starting from the first display, the same proof of part (a) of Lemma 1 applies?)}

In order to prove part (a) of the corollary, assume $0\leq x\leq x^{*}$ and set $s_k=\lfloor n-x\sqrt{n/2}\rfloor$. From (\ref{multi-sum}), we find
\begin{align*}
&p_n(\ell_k=s_k) = \\ 
&\sum_{s_1=1}^{s_k^{*}}\sum_{s_2=s_1}^{s_k^*}\cdots\sum_{s_{k-1}=s_{k-2}}^{s_k^{*}}\left[\prod_{\ell=0}^{k-2}\mu_{n-\ell}(s_{\ell+1}) p_{n-k+1}(\ell_1=s_k)+\sum_{\omega\neq \mu \mu \cdots \mu}\left(\prod_{\ell=0}^{k-2}\omega^{[\ell]}_{n-n_{\omega,\ell}-\ell}(s_{\ell+1}-n_{\omega,\ell})\right) p_{n-n_{\omega,k-1}-k+1}(\ell_1=s_k)\right]. 
\end{align*}
Then, the expansion of Lemma 1 for the factors $p_{n-n_{\omega,k-1}-k+1}(\ell_1=s_k)$ coupled with Lemmas \ref{ll-longest} and \ref{le4} yield 
\begin{align*}
&p_n(\ell_k=s_k) = \\
&\left[\frac{x}{\sqrt{n/2}}e^{-x^2/2}(1+o(1))+{\mathcal O}\left(\frac{e^{-x^2/2}}{n}\right) \right]\sum_{s_1=1}^{s_k^{*}}\sum_{s_2=s_1}^{s_k^*}\cdots\sum_{s_{k-1}=s_{k-2}}^{s_k^{*}}\left[\prod_{\ell=0}^{k-2}\mu_{n-\ell}(s_{\ell+1}) +\sum_{\omega\neq \mu \mu \cdots \mu}\left(\prod_{\ell=0}^{k-2}\omega^{[\ell]}_{n-n_{\omega,\ell}-\ell}(s_{\ell+1}-n_{\omega,\ell})\right) \right] \\
&=\left[\frac{x}{\sqrt{n/2}}e^{-x^2/2}(1+o(1))+{\mathcal O}\left(\frac{e^{-x^2/2}}{n}\right) \right] \left[  \frac{x^{2k-2}}{2^{k-1}(k-1)!}+{\mathcal O}\left(\frac{1+ x^{2k-3}}{\sqrt{n}}\right)  +  {\mathcal O}\left(\frac{1+ x^{2k-1}}{\sqrt{n}}\right) \right]  \\
&=\frac{x^{2k-1}}{2^{k-1}(k-1)!\sqrt{n/2}}e^{-x^2/2}(1+o(1))+{\mathcal O}\left(\frac{e^{-x^2/2}}{n}\right),
\end{align*}
as claimed in (a).
%{\bf Is this what you mean? Why this error term? $x$ could be $n^{1/7}$ with $k$ as large as we want.}

For part (b) we can write $p_n(\ell_{k}\leq n-x^{*}\sqrt{n/2}) = \sum_{x} p_n(\ell_{k} = \lfloor n-x\sqrt{n/2} \rfloor) = \sum_{x} p_n(\ell_{k} = s_k)$, where the sum proceeds in steps of $\sqrt{2/n}$ over the range $x^{*} \leq x \leq \sqrt{2n}$ and we set $s_k=\lfloor n-x\sqrt{n/2}\rfloor$.
Hence, by using (\ref{multi-sum}) together with Lemmas \ref{ll-longest} and \ref{le4}, we obtain
\begin{align*}
& p_n(\ell_{k}\leq n-x^{*}\sqrt{n/2}) =
\sum_x \sum_{s_1=1}^{s_k^{*}}\sum_{s_2=s_1}^{s_k^*}\cdots\sum_{s_{k-1}=s_{k-2}}^{s_k^{*}}\sum_{\omega}\left(\prod_{\ell=0}^{k-2}\omega^{[\ell]}_{n-n_{\omega,\ell}-\ell}(s_{\ell+1}-n_{\omega,\ell})\right)p_{n-n_{\omega,k-1}-k+1}(\ell_1=s_k) \\
&= \sum_{\omega} \sum_{x} \sum_{s_1=1}^{s_k^{*}}\sum_{s_2=s_1}^{s_k^*}\cdots\sum_{s_{k-1}=s_{k-2}}^{s_k^{*}}\left(\prod_{\ell=0}^{k-2}\omega^{[\ell]}_{n-n_{\omega,\ell}-\ell}(s_{\ell+1}-n_{\omega,\ell})\right)  p_{n-n_{\omega,k-1}-k+1}(\ell_1=s_k) \\
&= \sum_x \sum_{s_1=1}^{s_k^{*}}\sum_{s_2=s_1}^{s_k^*}\cdots\sum_{s_{k-1}=s_{k-2}}^{s_k^{*}}\left(\prod_{\ell=0}^{k-2}\mu_{n-\ell}(s_{\ell+1})\right)  p_{n-k+1}(\ell_1=s_k)\\
&+ \sum_{\omega \neq \mu \mu \cdots \mu} \sum_x \sum_{s_1=1}^{s_k^{*}}\sum_{s_2=s_1}^{s_k^*}\cdots\sum_{s_{k-1}=s_{k-2}}^{s_k^{*}}\left(\prod_{\ell=0}^{k-2}\omega^{[\ell]}_{n-n_{\omega,\ell}-\ell}(s_{\ell+1}-n_{\omega,\ell})\right)  p_{n-n_{\omega,k-1}-k+1}(\ell_1=s_k) \\
&= \sum_x p_{n-k+1}(\ell_1=s_k) \sum_{s_1=1}^{s_k^{*}}\sum_{s_2=s_1}^{s_k^*}\cdots\sum_{s_{k-1}=s_{k-2}}^{s_k^{*}}\left(\prod_{\ell=0}^{k-2}\mu_{n-\ell}(s_{\ell+1})\right)  \\
&+ \sum_{\omega \neq \mu \mu \cdots \mu} \sum_x p_{n-n_{\omega,k-1}-k+1}(\ell_1=s_k) \sum_{s_1=1}^{s_k^{*}}\sum_{s_2=s_1}^{s_k^*}\cdots\sum_{s_{k-1}=s_{k-2}}^{s_k^{*}}\left(\prod_{\ell=0}^{k-2}\omega^{[\ell]}_{n-n_{\omega,\ell}-\ell}(s_{\ell+1}-n_{\omega,\ell})\right)   \\
&=  \sum_x p_{n-k+1}(\ell_1=s_k) \left[   \frac{x^{2k-2}}{2^{k-1}(k-1)!}+{\mathcal O}\left(\frac{1+ x^{2k-3}}{\sqrt{n}}\right) \right] 
+ \sum_{\omega \neq \mu \mu \cdots \mu} \sum_x p_{n-n_{\omega,k-1}-k+1}(\ell_1=s_k) \left[   {\mathcal O}\left(\frac{1+ x^{2k-1}}{\sqrt{n}}\right) \right]. 
\end{align*}
Finally, since $x \leq \sqrt{2n}$, we have
\begin{align*}
p_n(\ell_{k}\leq n-x^{*}\sqrt{n/2}) &=  {\mathcal O}\left(n^{k-1}\right) \sum_x p_{n-k+1}(\ell_1=s_k) 
+  {\mathcal O}\left(n^{k-1}\right) \sum_{\omega \neq \mu \mu \cdots \mu} \sum_x p_{n-n_{\omega,k-1}-k+1}(\ell_1=s_k)  \\
&=  {\mathcal O}\left(n^{k-1}\right) p_{n-k+1}(\ell_1 \leq n^{*}) 
+  {\mathcal O}\left(n^{k-1}\right) \sum_{\omega \neq \mu \mu \cdots \mu} p_{n-n_{\omega,k-1}-k+1}(\ell_1 \leq n^{*})  \\
&=  {\mathcal O}\left(n^{k-1}\right) {\mathcal O}\left(e^{-n^{2/7}/2}\right) =  {\mathcal O}\left(n^{k-1} e^{-n^{2/7}/2}\right). \,\,\, \qed
\end{align*}

The next theorem, which extends Proposition~\ref{limit-law}, is our main result.
\begin{thm}\label{teo}
For a fixed $k \geq 1$, let $\ell_k$ be the $k$th largest external branch length in a random ordered history of size $n$ selected uniformly at random and denote by $\chi(2k)$ the $\chi$-distribution with $2k$ degrees of freedom. Then, as $n\rightarrow\infty$,
\[
\frac{n-\ell_k}{\sqrt{n/2}}\stackrel{d}{\longrightarrow}\chi(2k),
\]
with convergence of all moments. In particular, the mean and the variance of $\ell_k$ satisfy respectively
\begin{equation} \label{mean-var-gen}
{\mathbb E}(\ell_k)\sim n\qquad\text{and}\qquad{\rm Var}(\ell_k)\sim\left(k-\frac{\pi k^2}{16^k}\binom{2k}{k}^2\right)n.
\end{equation}
\end{thm}
\pf Following the proof of Proposition \ref{limit-law}, we show that all moments converge, which implies convergence in distribution. Starting from
\[
{\mathbb E}\left(\frac{n-\ell_k}{\sqrt{n/2}}\right)^m=\sum_{s=0}^{n}\left(\frac{n-s}{\sqrt{n/2}}\right)^m p_n(\ell_k=s),
\]
we replace $s$ by $s=n-x\sqrt{n/2}$ and break the sum into two parts obtaining
$$\sum_{x=0}^{\sqrt{2n}} x^m p_n(\ell_k=n-x\sqrt{n/2} )=\sum_{0 \leq x< n^{1/7}} x^m p_n(\ell_k=n-x\sqrt{n/2} ) +\sum_{n^{1/7} \leq x \leq \sqrt{2n}} x^m p_n(\ell_k=n-x\sqrt{n/2} ) \equiv \Sigma_1+\Sigma_2,$$
where all the sums proceed in steps of size $\sqrt{2/n}$. 
For $\Sigma_2$, by part (b) of the latter corollary, we have
\[
\Sigma_2={\mathcal O}\left(n^{m/2+k-1}e^{-n^{2/7}/2}\right)=o(1).
\]
For $\Sigma_1$, by part (a) of Corollary \ref{coro1}, we have
\[
\Sigma_1=\frac{1+o(1)}{2^{k-1} (k-1)!} \cdot \sum_{0\leq x<n^{1/7}}\frac{x^{m+2k-1}}{\sqrt{n/2}}e^{-x^2/2}+{\mathcal O}\left(n^{-1}\sum_{0\leq x<n^{1/7}}e^{-x^2/2}\right).
\]
Hence, the Riemann sum in $\Sigma_1$ can be approximated by the integral $\int_{0}^{n^{1/7}}x^{m+2k-1}e^{-x^2/2}{\rm d}x$, which converges to $\int_{0}^{\infty}x^{m+2k-1}e^{-x^2/2}{\rm d}x$. Overall,
\begin{align*}
{\mathbb E}\left(\frac{n-\ell_k}{\sqrt{n/2}}\right)^m & \stackrel{n \to \infty}{\longrightarrow} \frac{1}{2^{k-1} (k-1)!} \int_{0}^{\infty}x^{m+2k-1}e^{-x^2/2}{\rm d}x = \frac{1}{2^{k-1} (k-1)!} \cdot 2^{m/2+k-1} \, \Gamma\left( \frac{m+2k}{2} \right) \\
& = 2^{m/2} \, \frac{\Gamma\left( \frac{m}{2} + k\right)}{\Gamma(k)}
\end{align*}
which proves the claimed convergence of moments. Finally, (\ref{mean-var-gen}) follows from this convergence by straightforward computation. For instance, setting $m=1$ we obtain 
\begin{equation}\label{mean}
\frac{n-{\mathbb E}(\ell_k)}{\sqrt{n/2}} \stackrel{n \to \infty}{\longrightarrow} \frac{\sqrt{2 \pi} k {{2k}\choose{k}} }{4^k},
\end{equation}
and similarly for the variance.
\qed

%\noindent{\bf****************************************STOP****************************************}

\section{Conclusions}

\begin{figure}[tpb]
\begin{center}
\includegraphics*[scale=0.57,trim=0 0 0 0]{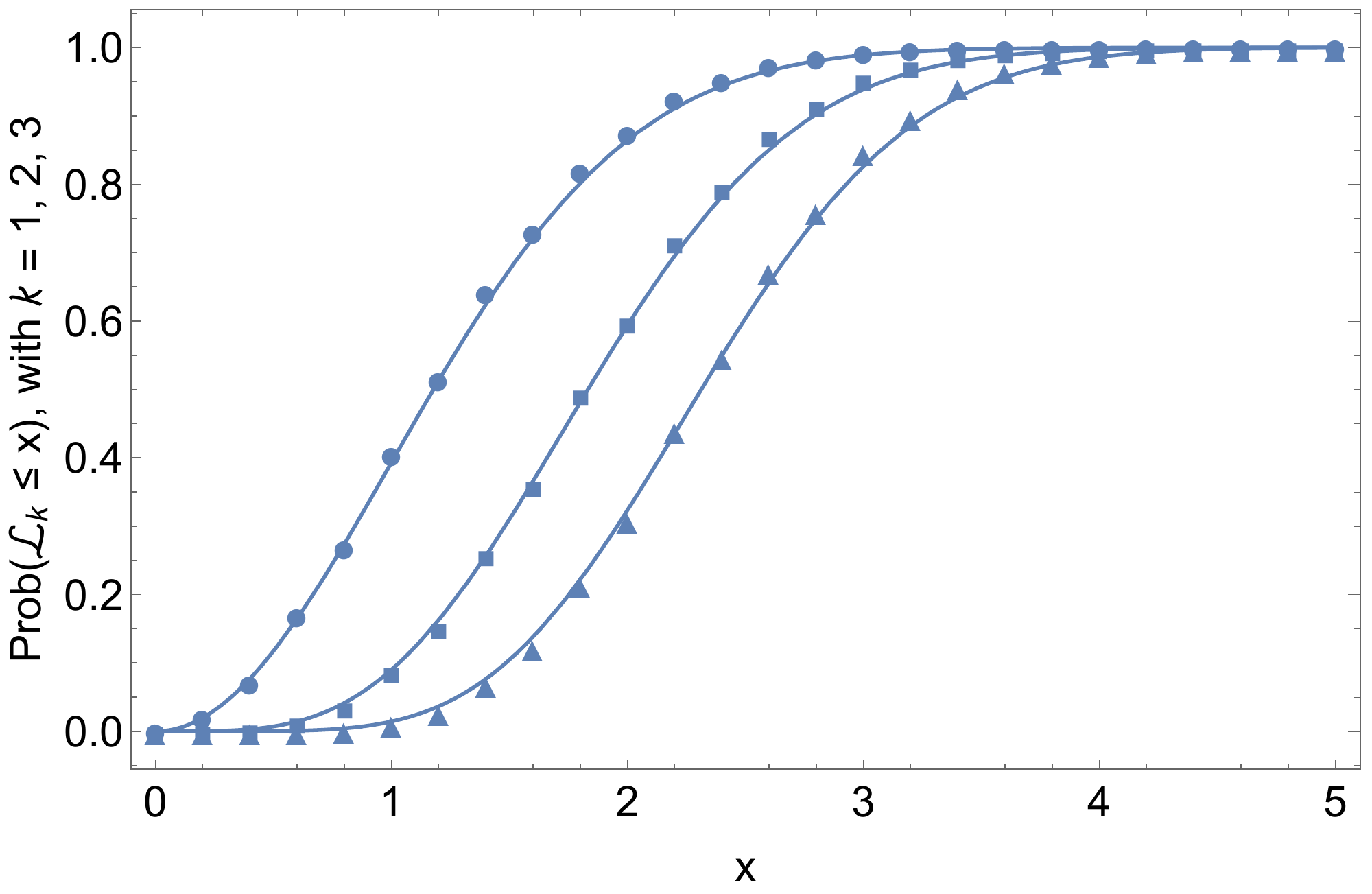}
\end{center}
\vspace{-.7cm}
\caption{{\small Probability that for $n=1000$ the rescaled variable $\mathcal{L}_k \equiv \frac{n-\ell_k}{\sqrt{n/2}}$ is less than or equal to $x \in [0,5]$ (in steps of $0.2$), when $k=1$ (dots), $k=2$ (squares), and $k=3$ (triangles). Values are calculated from Eqs. (\ref{prop1a}), (\ref{prop1b}), and (\ref{prop1c}). Solid lines give the cumulative function for the $\chi$-distribution with $2k$ degrees of freedom, with $k=1, 2, 3$ from left to right.  
}} \label{figcumul}
\end{figure} 

For random histories of fixed size $n$ selected under the Yule probability model, or, equivalently, for ordered histories of size $n$ selected uniformly at random, we have studied the variable $\ell_k$ defined as the $k$th largest length of an external branch. Measuring the length of an external branch as the rank of its parent node, Theorem \ref{teo} shows that the rescaled variable $\mathcal{L}_k \equiv \frac{n-\ell_k}{\sqrt{n/2}}$ follows asymptotically a $\chi$-distribution with $2k$ degrees of freedom (Fig. \ref{figcumul}), with convergence of all moments. The mean of $\ell_k$ is shown to be asymptotically equivalent to $n$, independently of $k$. More precisely, by plugging the approximation ${{2k}\choose{k}} \approx \frac{4^k}{\sqrt{\pi k}}$ into (\ref{mean}), we find that ${\mathbb E}(\ell_k)$ behaves like $n - \sqrt{k \, n}$ for increasing $n$. The variance of $\ell_k$ is asymptotically equivalent to $\left(k-\frac{\pi k^2}{16^k}\binom{2k}{k}^2\right)n$. 

Our approach has used a well known correspondence between trees and permutations, in which the $k$th largest length of an external branch of an ordered history of size $n$ is the value of the $k$th largest non-peak entry in the associated permutation of size $n-1$ (Section \ref{sec:2}). Thus, Proposition \ref{firstprop} and Theorem \ref{teo} also contribute to the study of the probabilistic properties of the value-peaks of permutations investigated in \cite{peaks}. 

In this paper we focused only on the {\it discrete} length of the external branches of random trees. Nevertheless, our results can also find applications in the analysis of the {\it time} length of the external branches of ``coalescent'' trees \cite{hudson, kingman:1982, tajima}. A coalescent tree of size $n$ is a pair consisting of a random Yule history $t$ of $n$ leaves and a sequence $( \tau_2, \dots, \tau_{n} )$ of independent exponentially distributed random variables assigning a time length to the different layers of $t$ (Fig. \ref{layers}). The variable $\tau_i$ gives the time length of the layer in which exactly $i$ branches of $t$ coexist, and its mean is $\mathbb{E}(\tau_i) = 1/\lambda_i$, with $\lambda_i = {{i}\choose{2}}$. Hence, the expected value of the time length of an external branch of $t$ of discrete length $s$ can be calculated as 
%is the sum of the expectations of the time length of the last $s$ layers of the history $t$, that is,
$\sum_{i=n+1-s}^{n} \mathbb{E}(\tau_i) = \frac{2}{n-s} - \frac{2}{n}.$ By using our finding that ${\mathbb E}(\ell_k) \approx n - \sqrt{k \, n}$, we thus see that, in a random coalescent tree of large size $n$, the mean of the $k$th time length of an external branch will behave roughly like $\frac{2}{\sqrt{k \, n}}.$ 
%Note that for instance, we can derive an expansion for the mean:
%\[
%E\left(\sum_{i=n+1-\ell_{k}}^{n}\tau_i\right)=2{\mathbb E}(1/(n-\ell_k))-2/n.
%\]
%For the mean on the right-hand side, we can again use Corollary 1 which yields
%\[
%{\mathbb E}(1/(n-\ell_k))\sim\frac{\sqrt{\pi}}{\sqrt{n}}\frac{1}{4^{k-1}}\binom{2k-2}{k-1}
%\]
%which for large $k$ is $\approx 1/\sqrt{kn}$ and thus,the above mean is $\approx 2/\sqrt{kn}$ as you claimed. (I am not so sure whether subtracting $2/n$ gives a good approximation; there might be a term of the same shape coming from the mean of $1/(n-\ell_k)$).

Yule and coalescent trees enable the simulation of the spread of mutations in a population under neutral evolution.
Singleton mutations---i.e. mutations affecting single individuals---can be modeled as random events occurring along the external branches of the tree. Doubleton mutations---which affect pairs of individuals---take place along those branches of the tree from which exactly two leaves descend. 
It would be of interest to extend the calculations of this article to investigate the length of this additional type of branches. 

\medskip
\noindent
{\footnotesize
{\bf Acknowledgments} Support to MF was provided by the MOST (Ministry of Science and Technology, Taiwan) grant MOST-111-2115-M-004-002-MY2.}


\begin{thebibliography}{5}

\bibitem{blum}
M.~G.~B. Blum, O. Fran\c cois, {\it Minimal clade size and external branch length under the neutral coalescent}, Adv. Appl. Probab. 37 (2005): 647--662.

\bibitem{peaks}
P. Bouchard, H. Chang, J. Ma, J. Yeh, Y.~N. Yeh, {\it Value-peaks of permutations}, Electron. J. Comb. 17 (2010): article \# R46.

\bibitem{caliebe}
A. Caliebe, R. Neininger, M. Krawczak, U. R\"{o}sler, {\it On the length distribution of external branches in coalescence trees: genetic diversity within species}, Theor. Popul. Biol. 72 (2007): 245--252.

\bibitem{dahmer}
I. Dahmer, G. Kersting, {\it The internal branch length of the Kingman coalescent}, Ann. Appl. Probab. 25 (2015): 1325--1348.

\bibitem{diehl}
C. Diehl, G. Kersting, {\it External branch lengths of $\Lambda$-coalescents without a dust component}, Electron. J. Probab. 24 (2019): 1--36.

\bibitem{anconfig}
F. Disanto, M. Fuchs, A.~R. Paningbatan, N.~A. Rosenberg, {\it The distributions under two species-tree models of the
number of root ancestral configurations for matching gene trees and species trees}, Ann. Appl. Probab. (to appear), arXiv:2006.09106 (preprint).

\bibitem{DisantoAndWiehe}
F. Disanto, T. Wiehe, {\it Measuring the external branches of a Kingman tree: A discrete approach}, Theor. Popul. Biol. 134 (2020): 92--105.

\bibitem{freund}
F. Freund, M.  M{\"o}hle, {\it On the time back to the most recent common ancestor and the external branch length of the Bolthausen-Sznitman coalescent}, Markov Process. Related Fields 15 (2009): 387--416.

\bibitem{fu}
Y.~X. Fu, {\it Statistical tests of neutrality of mutations}, Genetics 133 (1993): 693--709.

\bibitem{goulden}
I.~P. Goulden, D.~M. Jackson, {\it Combinatorial Enumeration}, Wiley, Chichester (1983).

\bibitem{harding} 
E.~F. Harding, {\it The probabilities of rooted tree-shapes generated by random bifurcation}, Adv. Appl. Probab. 3 (1971): 44--77.

\bibitem{hudson} 
R.~R. Hudson, {\it Gene genealogies and the coalescent process}, Oxf. Surv. Evol. Biol. 7 (1990): 1--44.

\bibitem{janson}
S. Janson, G. Kersting, {\it On the total external length of the Kingman coalescent}, Electron. J. Probab. 16 (2011): 2203--2218.

\bibitem{kingman:1982}
J.~F.~C.~Kingman,
{\it The coalescent},k
Stoch.~Proc.~Appl. 13 (1982): 235--248.

\bibitem{NielsenAndSlatkin}
R. Nielsen, M. Slatkin, {\it An Introduction to Population Genetics: Theory and Applications} Sinauer Associates, Sunderland, Massachusetts (2013).


\bibitem{rosenberg} 
N.~A. Rosenberg, {\it The mean and variance of the numbers of r-pronged nodes and r-caterpillars in Yule-generated genealogical trees}, Ann. Comb. 10 (2006): 129--146.

\bibitem{tajima} 
F. Tajima, {\it Evolutionary relationship of DNA sequences in finite populations}, Genetics 105 (1983): 437--460.

\bibitem{yule} 
G.~U. Yule, {\it A mathematical theory of evolution based on the conclusions of Dr. J.C. Willis, F.R.S.}, Philos. Trans. Roy. Soc. Lond. Ser. B 213 (1924): 21--87.


\end{thebibliography}
\end{document}